\documentclass[12pt,fleqn]{article}

\usepackage{amsmath}
\usepackage{amssymb}
\usepackage{amsfonts}
\usepackage{amsthm}
\usepackage{latexsym}
\usepackage{color}
\usepackage{graphicx}

\newtheorem{thm}{Theorem}[section]

\newtheorem{lem}{Lemma}[section]

\def\d{\, \mathrm d}
\def\Re{{\mathbb R}}

\newcommand{\Var}{\operatorname{Var}}

\renewcommand{\P}{\operatorname{P}}
\renewcommand{\mathbf}{\boldsymbol}
\newcommand{\iid}{\operatorname{iid}}

\newcommand{\ED}{{\mathbb F}_n}
\def\op{o_{\rm p}}

\def\Var{\mathop{\rm Var}\nolimits}
\def\Bias{\mathop{\rm Bias}\nolimits}

\renewcommand{\bf}{\bfseries}
\renewcommand{\tilde}[1]{\widetilde{#1}}
\renewcommand{\hat}[1]{\widehat{#1}}

\def\bea{\begin{eqnarray*}}
\def\eea{\end{eqnarray*}}
\def\be{\begin{equation}}
\def\ee{\end{equation}}
\def\bean{\begin{eqnarray}}
\def\eean{\end{eqnarray}}

\begin{document}

\title{\bf Smooth tail index estimation}
\author{Samuel M\"{u}ller$^{\text{a}}\footnote{Corresponding author. Email: mueller@maths.uwa.edu.au (S. M\"uller), tel: +61 (0)8 6488 7173}$, \ \ Kaspar Rufibach$^{\text{b}}\footnote{K.\ Rufibach (kaspar.rufibach@stanford.edu) is partially supported by Swiss National Science foundation. }$\\
$^{\text{a}}$School of Mathematics and Statistics, University of Western Australia\\ Crawley, WA 6009, Australia\\
$^{\text{b}}$Department of Statistics, Stanford University\\ Stanford, CA 94305, United States\\
}\date{\today}\maketitle

\abstract{Both parametric distribution functions appearing in extreme value theory - the generalized extreme
value distribution and the generalized Pareto distribution - have log-concave densities if the extreme value
index $\gamma\in[-1,0]$. Replacing the order statistics in tail index estimators by their corresponding
quantiles from the distribution function that is based on the estimated log-concave
density $\hat f_n$ leads to novel smooth quantile and tail index estimators. These new estimators aim at
estimating the tail index especially in small samples. Acting as a smoother of the empirical distribution function, the
log--concave distribution function estimator reduces estimation variability to a much greater extent than it
introduces bias. As a consequence, Monte Carlo simulations demonstrate
that the smoothed version of the estimators are well superior to their non-smoothed counterparts, in terms of
mean squared error.

\medskip


\medskip


\section{Introduction} \label{section: introduction}
It is a well--known fact that asymptotic results are in general at best approximately valid in small-sample
problems, but that in the latter situation bias is often a serious issue. For example in extreme value theory
the small-sample bias in the estimation of the tail index is severely. We refer to \cite{Dacorogna_1995} for a
study of a number of estimators. There are only a few more articles that focus on the small-sample performance
of tail-index estimators. We are aware of  \cite{Brazauskas_2001}, \cite{Huisman_2001} and \cite{Wagner_2004}.
These consider mainly tail-index estimation for heavy-tailed distributions. In this article we investigate the
small-sample behavior of a new smooth tail-index estimator for thin-tailed and distributions with finite
endpoint. The main aim of this article is to introduce the smoothed estimators that exploit the log-concavity of
the limiting density of the exceedances or of the largest order statistics, respectively. In Section \ref{Sec1}
we present the connection of log-concavity and extreme-value theory. Further, we show in Section \ref{section:
tail index estimation} that replacing the empirical distribution function by the smooth estimator $\hat F_n$
based on the log-concave density estimator - see (\ref{def: Fhat}) for a proper definition of $\hat F_n$ - leads
to novel tail-index estimators that exhibit substantially decreased mean-squared error in small-sample
situations. We illustrate this finding in a simulation study in Section \ref{section: simulations} for two
settings (1) a generalized Pareto distribution and (2) a domain of attraction scenario. The paper concludes with
some brief remarks in Section \ref{section: discussion}.

\section{Log-concavity in extreme-value theory}\label{Sec1}
\subsection{Max-domain of attraction of distributions with log-concave densities.}
 Let $\{X_n,\,n\geq 1\}$ be a sequence of independent random
variables with common distribution function $F$, such that $F$
belongs to the max-domain of attraction of $G$, denoted by $F\in\mathcal{D}(G)$, i.e.\ there
exist constants $a_n>0$, $b_n\in \Re$
such that for $x\in\Re$ and
\bea
G(x)&=&\lim_{n\rightarrow \infty} \P\Big(a^{-1}_n\cdot[\max(X_1,\ldots,X_n)+ b_n]\leq x\Big) \\
&=&\lim_{n\rightarrow \infty}F^{n}(a_n\cdot x-b_n)\nonumber
\eea which is equivalent to
\bea
 \sup_{x\in\Re} | F^{n}(a_n\cdot x-b_n) - G(x) |\rightarrow 0,\,\text{ as }n\rightarrow\infty.
\eea

\noindent From \cite{Gnedenko_1943} it is known that $F\in\mathcal{D}(G)$ if and only if
$G\in\{G_{\gamma}:\gamma\in\Re\}$, where
\begin{equation*}
G_{\gamma}(x)= \exp\left( -(1+\gamma x)^{-1/\gamma} \right), \label{eqn gev}
\;1+\gamma x >0,
\end{equation*}
where $G_{\gamma}$ is called the extreme value distribution with tail index $\gamma\in\Re$, shift parameter $0$, and scale parameter $1$. Since
\begin{equation*}
(1+\gamma x)^{-1/\gamma}\rightarrow\exp(-x), \;\text{ for } \gamma\rightarrow 0,
\end{equation*}
 interpret
$G_0(x)$ as $\exp(-e^{-x})$. The two most common settings in the analysis of extreme values are that we either have an observed sequence of independent and identically distributed maxima, $M_{n,1},\ldots,M_{n,k}$, or upper order statistics, $X_{(n)}\geq X_{(n-1)}\geq \ldots \geq X_{(n-K_{n})}$, from an independent and identically distributed sample $X_{1},\ldots,X_{n}$ with $K_{n} = k_{n} < n$ or $K_{n} = \#\{X_{i}: X_{i}\geq u_{n}\}$ a random number, where $u_n$ is some suitably chosen high threshold and $\#A$ denotes the number of elements in set $A$. We will focus on the latter setting and in general we assume that $K_{n}$ is an intermediate sequence that is $K_{n}/n \rightarrow 0$ and $K_{n}\rightarrow\infty$ if $n\rightarrow\infty$.

\subsection{The Generalized Pareto Distribution and log-concavity.}
To fix notation, define for a general distribution function $F$ the lower endpoint
$\alpha(F):=\inf\{x \in \Re :\;F(x)>0\}$ and the upper endpoint $\omega(F):=\sup\{x \in \Re:\;F(x)<1\}$. The quantile
function of $F$ for $q \in [0,1]$ is
\bea
F^{-1}(q):=\inf\{x \in \Re \: : \: F(x) \le q\}. \label{def: quantile fct}
\eea
Exceedances $X_{(n)}-u,\ldots,X_{(n-k_n+1)}-u$ of  a high threshold $u=u_{n}$ or of an intermediate order statistic $u=X_{(n-k_{n})}$ are typically modelled
by the generalized Pareto distribution (GPD), established by Pickands (see \cite{Pickands_1975}). For $\gamma \in \Re$ and
$\sigma \in (0,\infty)$ the density of the GPD is given by
\bea
    w_{\gamma,\sigma} (x) &:=& \sigma^{-1}(1+\gamma x/\sigma)^{-(1+1/\gamma)}, \ x \in[0,-\sigma/\gamma],
\eea where $w_{0,\sigma}$ is again defined via continuity: $w_{0,\sigma}(x) = \sigma^{-1} \exp(-x/\sigma)$ for $x \in [0, \infty)$.

\cite{Mueller_Rufi_2006a} showed that both parametric distribution functions appearing in extreme value theory -
the generalized extreme value distribution (GEV) and the generalized Pareto distribution - have log-concave densities
if the extreme value index $\gamma\in[-1,0]$ and that all distribution functions $F$ with log-concave density belong to
the max-domain of attraction of the GEV with $\gamma\in[-1,0]$. For any distribution function $F$ with corresponding
real-valued log-concave density function $f$, this latter $f$ can be written as
\bean
    f(x) &=& \exp \varphi(x), \label{def: logcon density}
\eean
for a concave function $\varphi \: : \: \Re \to [-\infty,\infty)$.
We will denote the class of distribution
functions $F$ having a log-concave density on its support $S=[\alpha(F),\omega(F)]$ by $\mathcal{F}_{\bar{\cap}}(S)$.

\subsection{The restriction $\gamma \in [-1,0]$ and connection to bump hunting.}
We are aware that the restriction of log-concave densities and of $\gamma \in [-1,0]$, respectively, is a drawback if
there is not sufficient evidence to assume that this restriction holds. On the other hand there are
often good reasons to assume that some distribution function $F$ has all its moments finite or that its support is
finite, implying $\gamma \in [-1,0]$.
Estimating the finite endpoint $\omega(F)$ of a distribution $F$ is linked to the problem of estimating $\gamma<0$ and its theory is well developed, see for example \cite{Hall_1982}, \cite{Csorgo_1989}, \cite{Hall_1999} and \cite{Ferreira_2003}.
The data set of the total life span of people who died in the Netherlands, who were born between the years 1877-1881, as analyzed in \cite{Aarssen_1994} is a real life
example that results in an estimated finite endpoint and an estimated tail index between $-1/2$ and $0$.
Further data sets on survival times of 208 mice exposed to radiation and on the men's 100m running times of the 1988 and 1992 Olympic Games are analyzed in \cite{Hall_1999}. By definition the distribution of the distance of two points in a closed
convex set has finite support and there are many open problems regarding its limit behavior.
For the current state of research we refer to \cite{Mayer_2006}. A further example that naturally leads to the restriction $\gamma<0$ is the estimation of the efficiency frontier in economics (see \cite{Farrell_1957}). In practical applications, $\gamma=-1/2$ is
often seen as natural lower bound, e.g.\ p.\ 62 in \cite{Kotz_2000}  or \cite{Ferreira_2003}.

In this context, we would like to point out \cite{rufibach_06_diss}, who proposes a multiscale procedure to identify collections
of intervals based on an i.i.d.\ sample where a density is either log-concave or log-convex. The chosen multiscale approach
ensures that the claims (i.e.~log-concavity or log-convexity) holds for all intervals in the collections simultaneously.
Additionally, it is shown that the method asymptotically keeps the level. This offers a way to ``pre-assess'' whether
$\gamma$ actually is in $[-1,0]$: For a chosen significance level $\alpha \in (0,1)$, apply this novel bump hunting method
to either the sample of exceedances or the whole sample of observations.
In the first case, if the collection of intervals claiming log-convexity is empty and the collection of intervals
claiming log-concavity at best contains an interval (almost) spanning the whole range of exceedances, then we could claim
with asymptotic probability $1-\alpha$ that indeed the observations we are looking at stem from a distribution
$F\in\mathcal{F}_{\bar{\cap}}\subset\mathcal{D}(G_\gamma)$, $\gamma \in [-1,0]$.

In the second case, we get an upper bound $m\geq k$ for the number of upper order statistics we should take into account.
Just define $m$ such that $X_{(n-m+1)}$ is the left-most endpoint of all intervals whereon the bump-hunting method claims
log-concavity of the underlying density.

In both cases, the new tail index estimation procedures as presented in Section \ref{section: tail index estimation} is adequate.

\section{Tail index estimation} \label{section: tail index estimation}
The estimation of $\gamma$ is besides the related high quantile estimation the most important problem in univariate extreme value theory and there exist a vast number of different approaches.
For example the Hill estimator (\cite{Hill_1975}),
the maximum likelihood estimator (\cite{Coles_1999}, \cite{Hall_1982}, \cite{Smith_1985a}, \cite{Smith_1987}, \cite{Smith_1985b}),
the moment estimator (\cite{Dekker_1989b}),
the (iterated) negative Hill estimator also known as Falk's estimator (\cite{Falk_1994b}, \cite{Falk_1995}, \cite{Mueller_2005}),
the (generalized) Pickands estimator (\cite{Pickands_1975}, \cite{Drees_1995}, \cite{Segers_2005}),
estimators based on near extremes (\cite{Mueller_2003}),
the weighted least squares estimator (\cite{Husler_2006}), probability weighted moments (\cite{Hosking_1985}),
and many more.
All these estimators are based on an intermediate sequence of upper order statistics and it is well known (see e.g.\
\cite{Groeneboom_2003}) that a major drawback of such estimators is their discrete character. Using
kernel-type estimators is one possibility to overcome this deficiency. We refer to \cite{Csorgo_1985}
for the smoothed Hill estimator in case $\gamma>0$ and to \cite{Groeneboom_2003} for general $\gamma\in\Re$. Our
alternative is to take advantage of the distribution function $\hat{F}_{n}$ based on the log-concave density estimator
$\hat f_n$, which is possible if the true $\gamma$ is in the restricted interval $[-1,0]$.

\subsection{Motivation of new estimators.} For an i.i.d.~sample $X_1,\ldots,X_n$ where $X_i$ has a log-concave
density function as introduced in (\ref{def: logcon density}), let $\ED$ be the empirical distribution function and
\bean
    \hat F_n (x)&:=& \int_{-\infty}^x \hat f_n (t) \d t \label{def: Fhat}
\eean be the smoothed distribution function based on the log-concave density estimator $\hat f_n = \exp \hat \varphi_n$.
A proper definition and
properties of $\hat f_n$ are given in \cite{rufibach_06_diss} and \cite{Rufibach_2004}. We only mention one special
feature: the estimator $  \exp \hat \varphi_n$ of the log-density is a piecewise linear function with knots only at some of the
observations points $X_1, \ldots, X_n$ and $\hat \varphi_n = - \infty$ on $\Re \setminus[X_1, X_n]$.
How to actually compute $\hat f_n$ is detailed in \cite{Rufibach_2006} and \cite{Duembgen_Huesler_Rufi_2006}.

In \cite{rufibach_06_diss} the following theorem is proven.

\begin{thm} \label{Theorem: log-concave distribution rate of convergence}
Let $X_1, \ldots, X_n$ be an i.i.d.~sample stemming from a distribution with log-concave density $f$ such that
 $f=\exp \varphi$ is H\"older-continuous for an exponent $\beta \in (1,2]$
and constant $L>0$. Furthermore, $\varphi'(x)-\varphi'(y) \ge C(y-x) \ $ for $ \ C>0 \ $ and $A\le x<y\le B$. Then, as
$n \to \infty$,
\bea
    \max_{t \in T_n} \, |\ED - \hat F_n|(t) \label{eq: conv2 dist}
    & = & \op \:  ( n^{-1/2} ),
\eea where
$T_n \to T:=[A,B]$. Furthermore, $\hat F_n(X_{(1)})=0$ and $\hat F_n (X_{(n)}) = 1$.
\end{thm}
This theorem implies that $\hat F_n$ is essentially equivalent
to $\ED$, but as the integral of a piecewise exponential function very smooth. These properties turn out to be
highly convenient in the estimation of the extreme value tail index $\gamma$. The smoothness of $\hat{F}_{n}$ reduces the variance not only considerably in the estimation of $\gamma$ but even for the estimation of quantiles of the generalized Pareto distribution, as is shown in Section  \ref{section: simulations}.

\medskip
Many well-known tail index estimators are based on a selection of log-spacings of the sample,
see also \cite{Segers_2005}.
The key idea is now simply to replace the order statistics (or quantiles of the empirical distribution function)
$X_{(i)}=\ED^{-1}(i/n), \ i=1,\ldots,n$ in these log-spacings by quantiles received via $\hat F_n$.
This yields modified versions of the uniformly minimum variance unbiased estimator
from Falk (\cite{Falk_1994b}, \cite{Falk_1995}) for the case of a known endpoint, the negative Hill estimator as defined in \cite{Falk_1995} and
Pickands' estimator (\cite{Pickands_1975}) for the case of a unknown endpoint. We will denote these new estimators
as ``smoothed estimators''.
We choose the first two estimators because of their outstanding performance for $\gamma<0$ and $\gamma<-1/2$,
respectively. On the other hand it is well known that Pickands' estimator is not efficient and in addition it is
able to estimate any $\gamma \in \Re$, thus the comparison of the original to the smoothed version is not entirely fair.
However, Pickands' estimator serves as the building block for much more efficient generalized Pickands' estimators that
are more general linear combinations of log-spacings of order statistics (\cite{Segers_2005}).

\subsection{Global and tail behavior.}Extreme value theory is, as the name suggests, tail focused. Hence, the behavior of the conditional distribution
$F_{X|X>u}$, where $u\rightarrow\omega(F)$, dominates the limit results. On the other hand, the log-concavity of
$f = F'$ is a strong assumption on the entire shape of the distribution function $F$. If this strong assumption holds
then the smoothing of the tail index estimators should be based on $\hat{F}_n$. Since tail index estimators use information of the upper tail of $\ED$ it was sufficient that the upper
tail of $F$ had log-concave density. Therefore, we investigate two settings in our simulation study in Section
\ref{section: simulations}. (1) $X_1,\ldots,X_n \iid \mathcal{L}(X)\in\mathcal{F}_{\bar{\cap}}(S)$. Here the data is
sharpened by $\hat{F}_n$. (2) $X_1,\ldots,X_n \iid \mathcal{L}(X)\in\mathcal{D}(G_\gamma; \gamma \in [-1,0])$. In this
more general situation $F$ is not necessarily in $\mathcal{F}_{\bar{\cap}}(S)$ . Before the smooth tail index estimators are computed, the range of log-concavity of $f$ has to be determined; either using bump-hunting or by imposing some assumptions on $F$. Here the data is sharpened by $\hat{F}_m$, where the smoothed distribution function is based on the $m>K_n$ largest order statistics and $m$ is such that $F_{X|X>X_{n-m}}\in \mathcal{F}_{\bar{\cap}}(S\cap[X_{n-m},\infty))$.

\subsection{Smooth tail index estimators.}
First, let us pin down some notation. Suppose we are given a sample $X_1, \ldots , X_n$ from a GPD wherefrom we
know that $\gamma \in [-1,0]$ and with empirical distribution function $\ED$. Denote the order statistics
by $X_{(1)}, \ldots, X_{(n)}$. For such a fixed sample,
define 
for $k = 4, \ldots, n$ and $H \in \{\ED, \hat F_n\}$:
\bea
    \hat \gamma^k_{\text{Pick}}(H) &=& \frac{1}{\log 2} \log \Bigl(\frac{H^{-1}((n-r_k(H)+1)/n)-H^{-1}((n-2r_k(H) +1)/n)}{H^{-1}((n-2r_k(H) +1)/n)-H^{-1}((n-4r_k(H) +1)/n)} \Bigr),
\eea
where
\bea
    r_k(H) &=&
    \begin{cases}
        \lfloor k/4 \rfloor & \text{if} \ \ H = \ED, \\
        k/4                 & \text{if} \ \ H = \hat F_n,
\end{cases}
\eea with $\lfloor m \rfloor := \max\{n \in \mathbb{N}_0 \: : \: n\le m\}$. This construction not only exploits the
superiority of the quantile estimates based on the smooth function $\hat F_n$, 
but also avoids ``rounding bias". Using the inverse of a continuous distribution function, quantiles do not
coincide for four consecutive $k$'s (order statistics), as it is the case for Pickands' original estimate.

To generalize the estimators in \cite{Falk_1994b} and \cite{Falk_1995}, no discrimination regarding continuity of
$H$ is necessary. For $H \in \{\ED, \hat F_n \}$ let
\bea
    \hat \gamma^k_{\text{Falk}}(H) &=& \frac{1}{k-1} \sum_{j=2}^k \log \Bigl(\frac{X_{(n)}-H^{-1}((n-j+1)/n)}{X_{(n)}-H^{-1}((n-k)/n)} \Bigr), \:\: k=3,\ldots,n-1 \\
    \hat \gamma^k_{\text{MVUE}}(H) &=& \frac{1}{k} \sum_{j=1}^k \log \Bigl(\frac{\omega(F)-H^{-1}((n-j+1)/n)}{\omega(F)-H^{-1}((n-k)/n)} \Bigr), \:\: k=2,\ldots,n-1
\eea where $F$ is the true distribution function of the $X_i$'s.
Note that $\hat \gamma^k_{\text{MVUE}}(H)$ is only consistent if $\gamma \in [-1, 0)$.
The chosen terminology reminds of the fact that when choosing $H=\ED$, the above estimators boil down to
Pickands', Falk's, and Falk's MVUE estimator as discussed at the beginning of this section.

The new, smooth tail index estimators are now simply $\hat \gamma^k_{\text{Pick}}(\hat F_n), \hat \gamma^k_{\text{Falk}}(\hat F_n)$, and
$\hat \gamma^k_{\text{MVUE}}(\hat F_n)$. Figure \ref{figure: Hill plots} displays Hill plots for two GPD pseudo-random samples, i.e.\ plots of the estimators
versus the number of order statistics $k$, for the smoothed and
unsmoothed versions for $n = 64$ and $\gamma \in \{-0.1, -0.75\}$.
\begin{figure}[!h]
\begin{center}
\vspace{-2cm}
\centerline{\includegraphics[width=16cm]{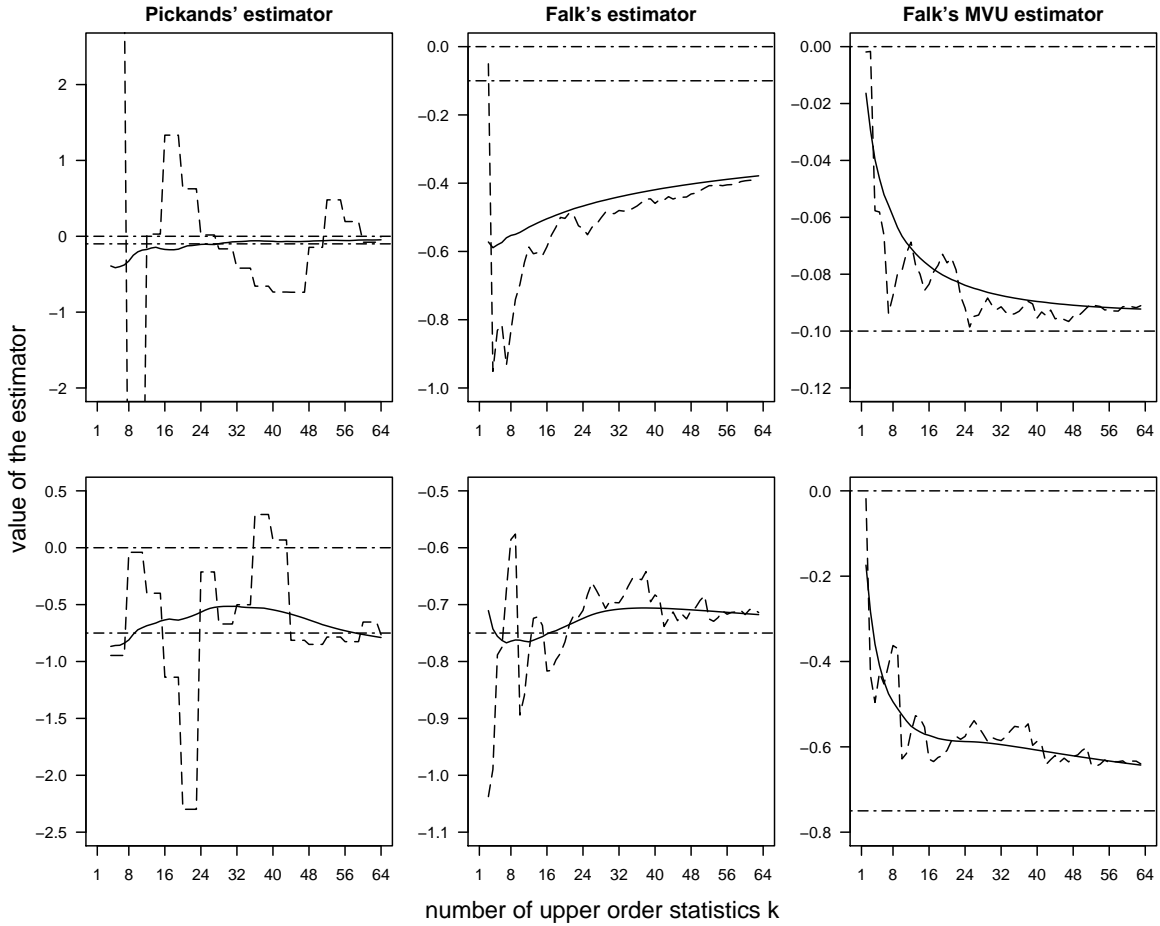}}
\end{center}
\vspace{-1cm}
\caption{Hill plots for $n = 64$ and $\gamma = -0.1$ (plots in the upper row) and $\gamma = -0.75$ (plots in the lower row), smoothed (--) and original (- -) versions.}
\label{figure: Hill plots}
\end{figure}

The smoothed estimators behave much more stable as
a function of $k$ and it is especially noteworthy that  for the two generated data sets all three smoothed estimators $\hat \gamma^k_{\text{Pick}}(\hat F_n)$,
$\hat \gamma^k_{\text{Falk}}(\hat F_n)$, $\hat \gamma^k_{\text{MVUE}}(\hat F_n) \in [-1,0]$ for every $k$.
By construction, $\hat{\gamma}_{\text{Pick}}:\Re^n\rightarrow \Re$ and
$\hat{\gamma}_{\text{Falk}},\hat{\gamma}_{\text{MVUE}}:\Re^n\rightarrow [-\infty, 0)$ which means that non permissible
estimates outside the interval $[-1,0]$ potentially occur. However, due to consistency of $\ED$ and $\hat F_n$,
this is asymptotically negligible. If in practice $\hat{\gamma}_{\cdot}\not\in [-1,0]$ then, a truncation to its closest
boundary value is recommendable if it is known that $\gamma\in[-1,0]$.

It is beyond the scope and not the primary goal of this article to discuss the asymptotic behavior of the
smoothed tail index estimators.


\pagebreak

\subsection{Further shape constraints.} Straightforward computation yields the following lemma.
\begin{lem}\label{lem: form of w_gam}
The density $w_{\gamma,\sigma}$ has the following qualitative properties which do not depend on the value of the
scale parameter $\sigma$:
\begin{table}[!h]
\begin{center}
\vspace{-0.3cm}
\begin{tabular}{ll}
\textit{property} & \textit{parameter range} \\ \hline
\textit{convex non-decreasing}  & $\gamma \le -1$ \\
\textit{concave non-increasing} & $\gamma \in [-1,-1/2]$ \\
\textit{convex non-increasing}  & $\gamma \ge -1/2$ \\ \hline
\textit{log-concave}            & $\gamma \in [-1,0]$ \\
\textit{log-convex}             & $\gamma \in (-\infty,-1]\cup[0,\infty)$
\end{tabular}
\end{center}
\caption{Form of $w_{\gamma, \sigma}$ as determined by the tail index $\gamma$.}
\label{table: further constraints}
\end{table}
\end{lem}

Since e.g. a density estimator $\tilde f_n$ in the class of convex decreasing densities is available, see \cite{Groeneboom_2001},
the latter lemma raises the possibility to define smooth estimators of the tail index for other ranges of $\gamma$.
However, especially in the latter case, slight modifications may be necessary, since
$\tilde F_n(X_{(n)}) \ne 1$. Furthermore, maximum likelihood estimators under other constraints may not be
as smooth as $\hat F_n$, since e.g. the estimator of a convex decreasing density is piecewise linear, whereas
for log--concave densities this form appears for the estimator of the log--density.

\subsection{Computational details.} These new smoothed estimators are made available in the R-package \verb"smoothtail" (see \cite{smoothtail}).
This latter package depends on the package \verb"logcondens" (\cite{logcondens}), which offers two algorithms for the (weighted) estimation
of an arbitrary log-concave density from an i.i.d.\ sample of observations. Both these packages are available from
CRAN.

\section{Simulations} \label{section: simulations}

\subsection{Estimation of quantiles.}
The computation of the non-smoothed estimators $\hat \gamma^k_{\text{Falk}}(\ED)$, $\hat \gamma^k_{\text{Pick}}(\ED)$ and
$\hat \gamma^k_{\text{MVUE}}(\ED)$ heavily relies on the order statistics  $X_{(i)}, i=1,\ldots,n$.
But, these simply estimate the quantiles $W_{\gamma, \sigma}^{-1}(i/n)$ of the distribution whereof we want to
estimate the tail index. Therefore, the accuracy of the tail index estimators is closely connected to the ability of
estimating these quantiles $W_{\gamma, \sigma}^{-1}(i/n)$.
To illustrate the superiority of log-concave quantile estimation
over simply taking order statistics, we calculated the relative efficiency of these two estimators. Since results were
similar over an extended range of $\sigma$'s, we concentrate on the case $\sigma = 1$.

To fix notation, define the log-concave estimate $\hat X_{(i)} = \hat F_{n}^{-1}(i/n)$ of an order statistic $X_{(i)}$,
for $i=1, \ldots, n$. Let $q_{(i)}$ denote either $\hat X_{(i)}$ or $X_{(i)}$, then its estimated variance and bias with
respect to the $i/n$-quantile of a GPD $W_{\gamma,1}$ for a fixed $\gamma \in [-1,0]$ given simulated $q_{(i),j}$ that are
based on $M$ generated samples
$X_{(1),j},\ldots,X_{(n),j}$, $j=1,\ldots,M$ of size $n$ drawn from $W_{\gamma,1}$ is defined as follows:
\bea
    \widehat{\Var}(q_{(i)},M)&:=& (M-1)^{-1} \sum_{j=1}^M \Big(q_{(i),j}- (1/M)\sum_{j=1}^M q_{(i),j}\Big)^2 \\
    \widehat{\Bias}(q_{(i)},\gamma,M)&:=& M^{-1} \Big(\sum_{j=1}^M q_{(i),j}\Big)-W_{\gamma,1}^{-1} (i/n).
\eea
The relative efficiency $\rho_{\gamma,n,M}(k)$ of log-concave quantile estimation to quantile estimation based on
order statistics is then
\bea
    \rho_{\gamma,n,M}(k)&=&\frac{[\widehat{\Bias}(\hat X_{(k)},\gamma,M)]^2+\widehat{\Var}(\hat X_{(k)},M)}{[\widehat{\Bias}(X_{(k)},\gamma,M)]^2+\widehat{\Var}(X_{(k)},M)}.
\eea

Figure \ref{figure: relative efficiencies quantiles} details $\rho_{\gamma,32,1000}(k)$ for
$\gamma \in \{-1, -0.75, -0.5, -0.25, 0\}$ as a function of $k$.
Relative efficiencies smaller than 1 are in favor of the log-concave quantile estimation
and indicate its superiority.
The use of the log-concave density estimator for the estimation of quantiles substantially reduces the variance of the
estimation, due to its smoothing property detailed in Theorem
\ref{Theorem: log-concave distribution rate of convergence}. This transfers to a reduced MSE, uniformly in
$\gamma$ and $k$, as is detailed in Figure \ref{figure: relative efficiencies quantiles}.

\begin{figure}[!h]
\begin{center}
\centerline{\includegraphics[width=14cm]{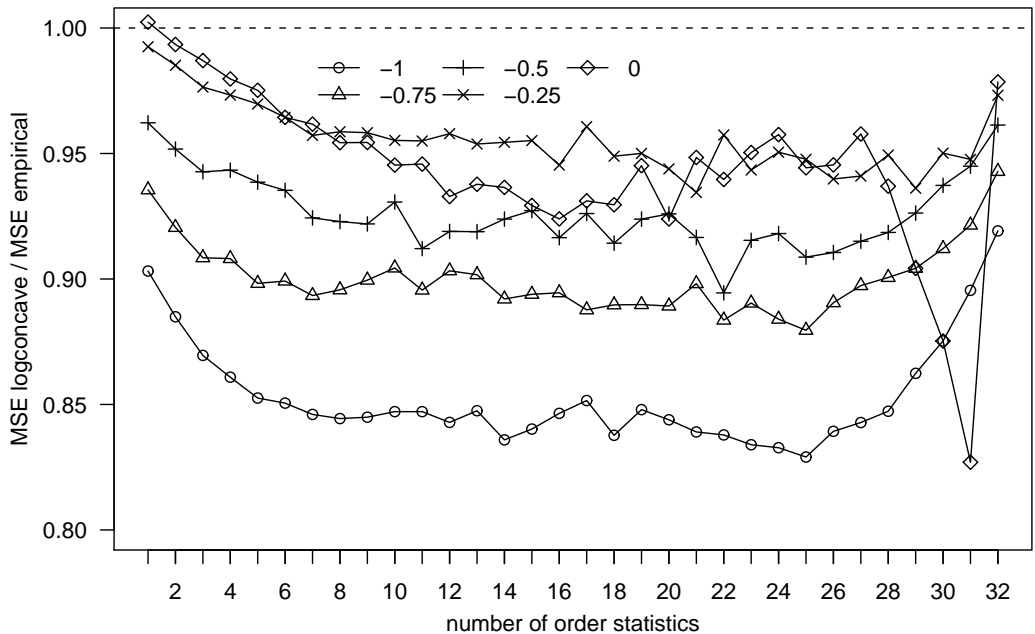}}
\end{center}
\vspace{-1.5cm}
\caption{\small Relative MSE in quantile estimation for $n = 32$.}
\label{figure: relative efficiencies quantiles}
\end{figure}

\subsection{Smoothed versus unsmoothed tail index estimators.}
To assess the effect of smoothing the tail index estimators, we perform a simulation study for two settings.
In Setting 1, we draw $M=1000$ samples of size $n_1 = 64$ from a GPD with $\sigma = 1$ and the extreme value tail index
$\gamma \in\{-1,-0.75,-0.5,-1/3,-0.25,-0.1\}$. For every $k$, the log-concave density is estimated based on the
full sample $X_{(1)}, \ldots, X_{(n)}$.
Setting 2 consists of $M=1000$ samples from a $\beta(\theta_1,\theta_2)$-distribution having density
\begin{equation}
f_{\theta_1, \theta_2}(x) = \frac{\Gamma(\theta_1+\theta_2)}{\Gamma(\theta_1)\Gamma(\theta_2)} x^{\theta_1 -1} (1-x)^{\theta_2 -1}\, , \ \ x \in [0,1]\, , \ \theta_1,\, \theta_2 > 0. \label{eqn: beta density 1}
\end{equation}
The upper tail of the $\beta$-distribution is dominated by $\theta_2$, since for $x\uparrow 1$ we have that
$f_{\theta_1, \theta_2}(x) = c(\theta_1,\theta_2)\cdot w_{-\theta_2^{-1},\theta_2^{-1}}(x)[1+o(1)]$. Thus,
$\beta(\theta_1,\theta_2) \in \mathcal{D}(G_{-1/\theta_2})$. We fix $\theta_1 = 1/2$ and thus only the upper tail of $f_{0.5, \theta_2}$ is log-concave.
For the simulations, we choose $\theta_2 = -\gamma^{-1} \in\{1,4/3,2,3,4,10\}$ and $n_2 = 128$ which
equals deliberately $2n_1$ in order to underline the difference between the two settings. Here the log-concave density
estimator is based on the largest $n_2/2 = 64$ order statistics. In both settings we present results for a single sample
size only because the results for $n_1, n_2/2 \in \{32,128,256,512\}$ are very similar.

Setting 1 represents the ``ideal'' of observing pure ``peak over threshold'' data
$$X_1,\ldots,X_n \iid \mathcal{L}(X) =\text{GPD}\in\mathcal{F}_{\bar{\cap}}(S).$$ Alternatively, Setting 2 stands for the more general
situation $$X_1,\ldots,X_n \iid \mathcal{L}(X)\in\mathcal{D}(G_\gamma; \gamma \in [-1,0]).$$
In the latter case, the well known problem of the tradeoff between bias and variance dominates the optimal choice of $k$ whereas in Setting 1,
$k$ can be chosen as large as possible since we are considering the ``perfect'' model.

\subsection{Simulation results.}
We compute relative efficiencies for the estimation of the tail index as for quantile estimation. For Setting 1, results
are displayed in Figure \ref{figure: relative efficiencies tail index n=64}. If one knows that $\gamma \in [-1,0]$,
then clearly using the smoothed estimator is most worthwhile for Pickands' estimator. However, also the other two
estimators are substantially improved, with highest gain in terms of MSE for small $k$'s.

\begin{figure}[!h]
\begin{center}
\vspace{-2.5cm}
\centerline{\includegraphics[width=16cm]{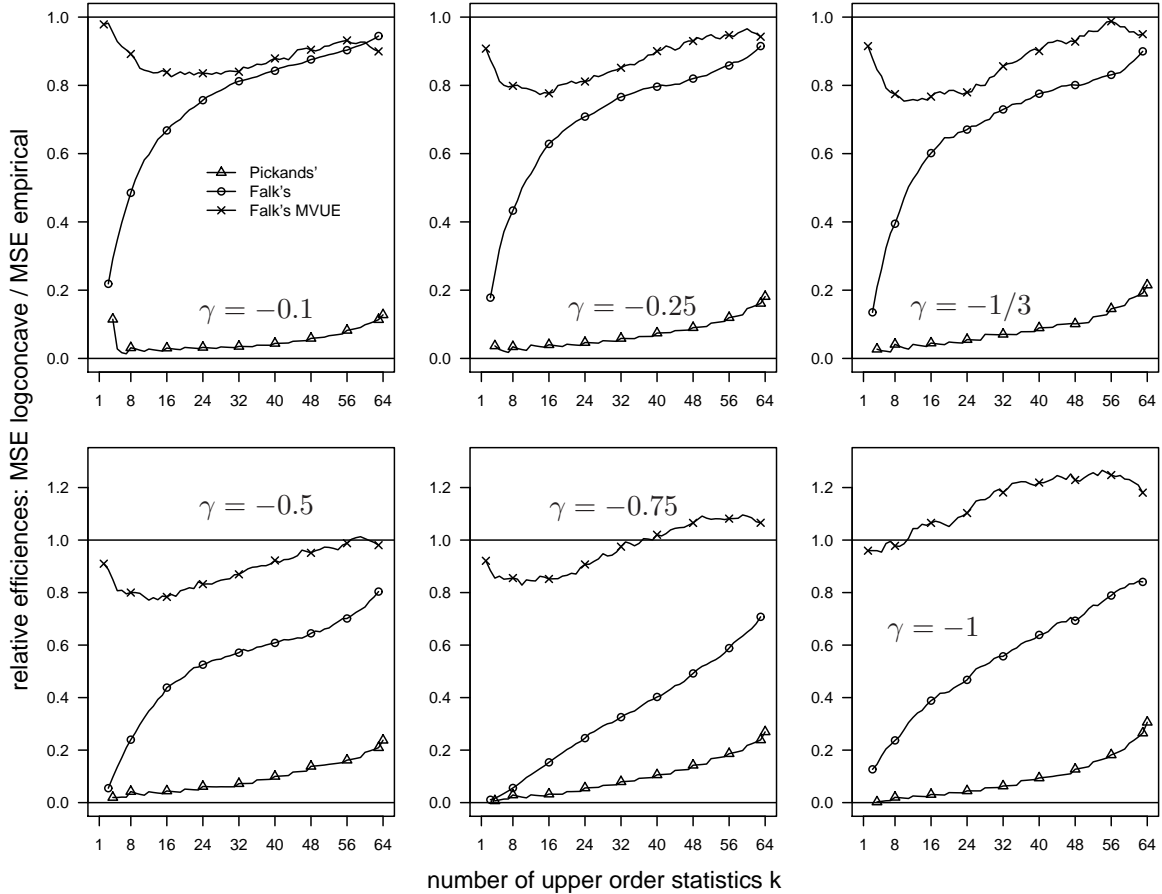}}
\vspace*{-8.5cm}\hspace{-9cm}{\small $\gamma = -0.1$}

\vspace*{2cm}\hspace{-9cm}{\small $\gamma = -0.5$}

\vspace*{-0.6cm}\hspace{0.5cm}{\small $\gamma = -0.75$}

\vspace*{1.0cm}\hspace{9cm}{\small $\gamma = -1$}

\vspace*{-4.9cm}\hspace{10cm}{\small $\gamma = -1/3$}

\vspace*{-0.6cm}\hspace{1cm}{\small $\gamma = -0.25$}

\vspace*{8.2cm}
\end{center}
\vspace{-1cm}
\caption{Relative MSE for tail index estimation for peaks over threshold data, $n = 64$.}
\label{figure: relative efficiencies tail index n=64}
\end{figure}

Figure \ref{figure: var bias and mse for setting 2} sheds light on the bias-variance trade-off in Setting 2 for
$\theta_2 = 3$; results for other choices of $\theta_2$ were absolutely similar and therefore omitted. The variance
in estimation of the tail index is dramatically reduced for the smoothed Pickands' estimate. The plot of
the bias against $k$ confirms that all estimators are biased as expected, since the data is generated by
(\ref{eqn: beta density 1}) and not by the GPD with $\gamma\in[-1,0]$. Especially for Pickands' estimator, the
bias is a much smoother function of $k$ for the smoothed estimators than this is the case for the original estimator.
\begin{figure}[!h]
\begin{center}
\centerline{\includegraphics[width=16cm]{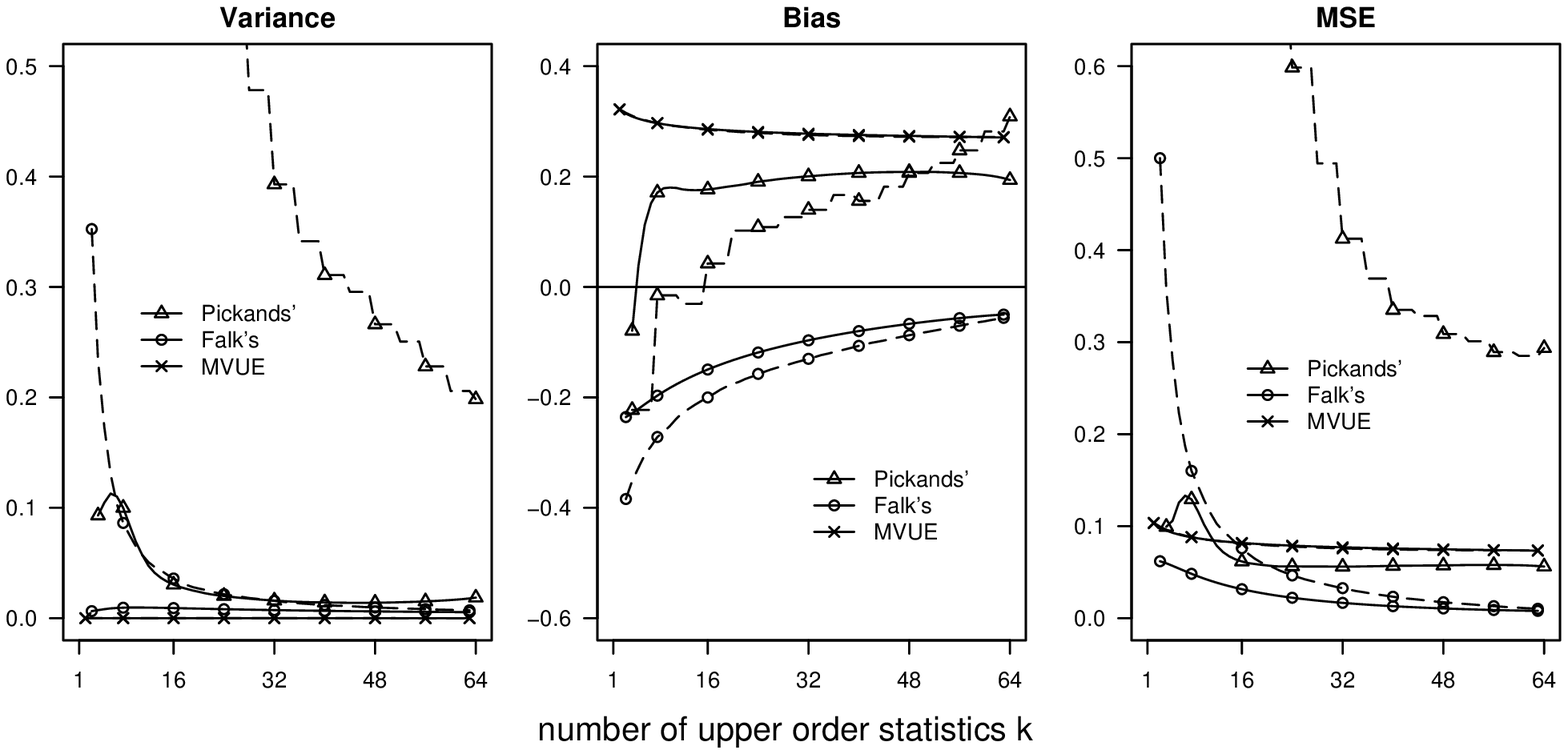}}
\end{center}
\vspace{-1.5cm}
\caption{\small Domain of attraction scenario: $n = 128$, $\theta_1 = 1/2$, $\theta_2 = 3$, smoothed (--) and original (- -) versions.}
\label{figure: var bias and mse for setting 2}
\end{figure}
Figure \ref{figure: relative efficiencies tail index beta n=128} shows the computed relative efficiencies for the
estimation of the tail index for Setting 2. The results for $\hat{\gamma}_{\text{Pick}}$ and $\hat{\gamma}_{\text{Falk}}$
are similar to those in Setting 1, yielding the most substantial improvement for Pickands' estimator. On the other hand the efficiency of the smoothed and original
$\hat{\gamma}_{\text{MVUE}}$ is almost 1, independent of $k$ and the choice of $\theta_2$.
\begin{figure}[h]
\begin{center}
\centerline{\includegraphics[width=16cm]{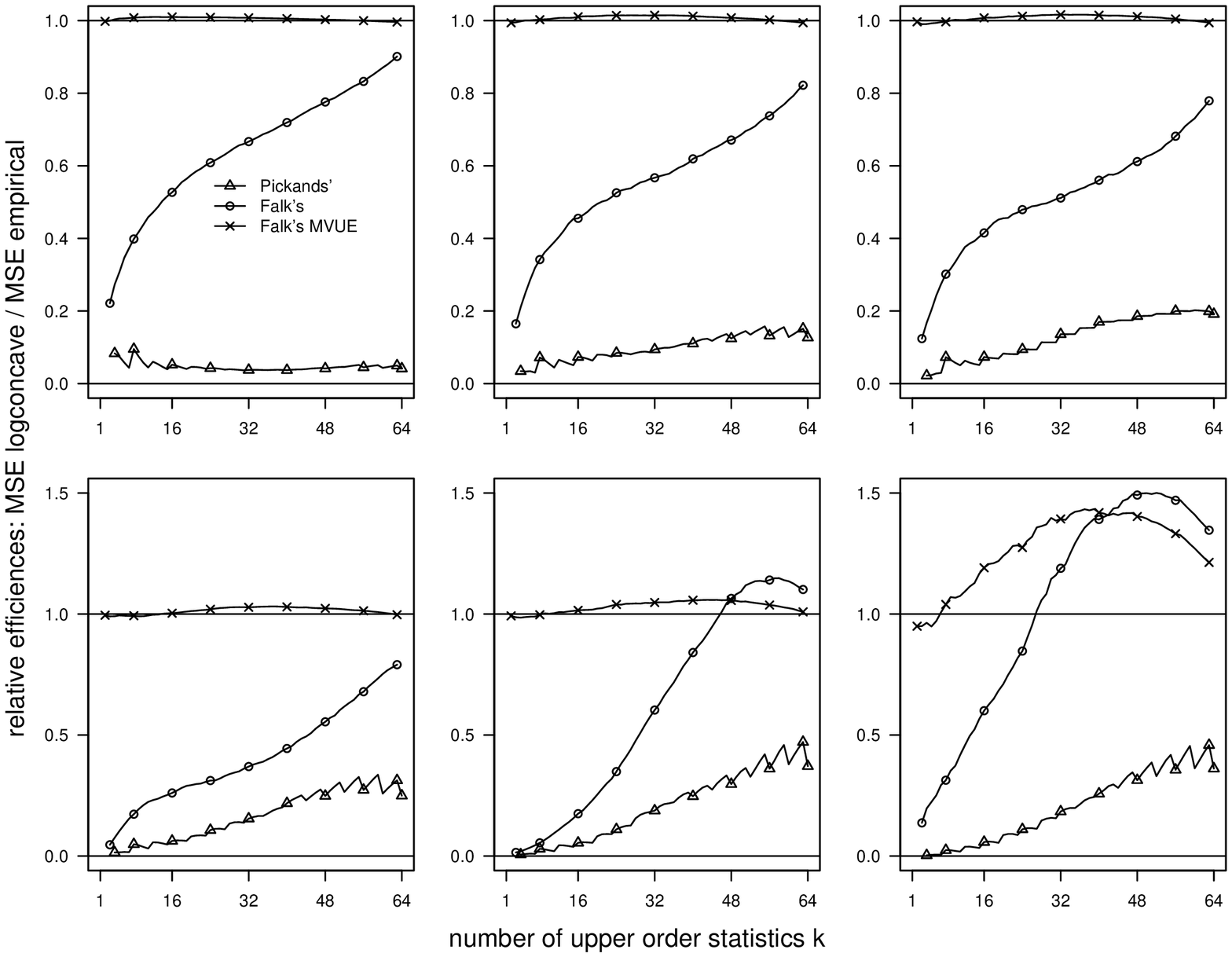}}
\end{center}
\vspace{-1,5cm}
\caption{\small Relative MSE in domain of attraction scenario, $n = 128$. $\theta_2$ equals 10, 4 and 3 in the top row and 2, $4/3$ and 1 in the bottom row.}
\label{figure: relative efficiencies tail index beta n=128}
\end{figure}

\section{Conclusions} \label{section: discussion}
In this article we showed that for the class of distributions $F\in\mathcal{D}(G_{\gamma}$) having log-concave
densities or log-concave conditional densities given that the observations exceed some threshold $u$,
respectively, the smoothing of the empirical distribution function by the corresponding log-concave density
estimator leads to improved quantile estimation uniformly in $\gamma$ and the quantile to be estimated and to
more efficient tail index estimators if it is known that $\gamma\in[-1,0]$. The reduction of the mean squared
error of the smoothed tail index estimator is not surprisingly most substantial for Pickands' estimator, given
the known poor efficiency of the latter. The estimator $\hat{F}_n$ smoothes the empirical distribution function
based on the global property of the log-concavity of the density of the underlying distribution function $F$. We
showed that if such a global property is present then it is of great value for the estimation of a tail property
such as the extreme value index $\gamma$. Of course the price to be paid is assuming that $\gamma$ is restricted
to the narrow interval $[-1,0]$ and it makes only sense to use the presented smoothed estimators if there is
sufficient reason for that assumption. However, the novel bump hunting method in \cite{rufibach_06_diss} is a
possible tool for testing if a log-concavity assumption holds which automatically causes the desired restriction
on $\gamma$. We want to stress the fact that the Hill plots in Section \ref{section: tail index estimation} and
the results in our simulation study are not based on an optimal choice of the setting parameters and are
reproducible for a wide range of parameter values for the $\beta(\theta_1,\theta_2)$-distribution unless the
upper tail of the density is log-concave. Our simulation study can be reproduced using the R-package
\verb"smoothtail" (see \cite{smoothtail}).

\subsection{Further research.} The results in this article raise many new open questions and problems. We only mention
two. At first it is challenging to prove asymptotic normality for the smoothed estimators. Then, by fitting the
generalized Pareto density to $\hat{f}_n$ or $\log w_{\gamma,\sigma}$ to $\hat{\varphi}_n$, respectively, opens the door
to construct novel estimators for the tail index that are worthwhile to be investigated in the future.

\section{Acknowledgments.} This research was initiated while both authors were working at the Institute of Mathematical
Statistics and Actuarial Science at the University of Bern.
K.\ Rufibach was partially supported by Swiss National Science foundation.

\end{document}